\documentclass{amsart}
\usepackage{amssymb}
\usepackage{a4}
\begin{document}
\newtheorem{theorem}{Theorem}[section]
\newtheorem{lemma}[theorem]{Lemma}
\newtheorem{example}[theorem]{Example}
\newtheorem{remark}[theorem]{Remark}
\def\ffrac#1#2{{\textstyle\frac{#1}{#2}}}
\def\qedbox{\hbox{$\rlap{$\sqcap$}\sqcup$}}
\makeatletter
 \renewcommand{\theequation}{%
 \thesection.\alph{equation}}
 \@addtoreset{equation}{section}
 \makeatother
\title[$k$-Curvature homogeneous manifolds]
{Complete $k$-curvature homogeneous pseudo-Riemannian manifolds}
\author{P. Gilkey and S. Nik\v cevi\'c}
\begin{address}{PG: Mathematics Department, University of Oregon,
Eugene Or 97403 USA.\newline Email: {\it gilkey@darkwing.uoregon.edu}}
\end{address}
\begin{address}{SN: Mathematical Institute, Sanu,
Knez Mihailova 35, p.p. 367,
11001 Belgrade,
Serbia and Montenegro.
\newline Email: {\it stanan@mi.sanu.ac.yu}}\end{address}
\begin{abstract} For $k\ge2$, we exhibit complete $k$-curvature homogeneous neutral signature
pseudo-Riemannian manifolds which are not locally affine homogeneous (and hence not locally homogeneous). All the local scalar
Weyl invariants of these manifolds vanish. These manifolds are Ricci flat, Osserman, and Ivanov-Petrova.\end{abstract}
\maketitle
\section{Introduction}\label{sect-1}

We consider a pair $\mathcal{M}:=(M,g_M)$ where $g_M$ is a pseudo-Riemannian metric of signature $(p,q)$ on a smooth manifold
$M$ of dimension $m:=p+q$. Let $R$ be the associated Riemann curvature tensor and let $\nabla^kR$ denote the $k^{\text{th}}$
covariant derivative of the curvature tensor. We say that $\mathcal{M}$ is $k$-curvature homogeneous if given any two points
$P,Q\in M$, there exists an isomorphism $\phi_{P,Q}$ from $T_PM$ to $T_QM$ so that
$$\phi^*g_Q=g_P,\ \phi^*R_Q=R_P,\ ...,\ \phi^*\nabla^kR_Q=\nabla^kR_P\,.$$
This means that the metric, curvature tensor, and covariant derivatives up to order $k$ of the curvature tensor ``look the
same'' at each point.

There is an equivalent algebraic formalism. Consider
$$\mathcal{U}_m^k:=(V,g,A^0,A^1,...,A^k)$$
where $g$ is an inner product on a $m$ dimensional real vector space $V$ and where $A^i\in\otimes^{4+i}V^*$ for $0\le i\le k$.
We say that
$\mathcal{U}_m^k$ is a {\it $k$-model} for $\mathcal{M}$ if for every point $P\in M$, there is an isomorphism
$\phi_P:T_PM\rightarrow V$ so that
$$\phi_P^*g=g_P,\ \phi_P^*A^0=R_P,\ ...,\ \phi_P^*A^k=\nabla^kR_P\,.$$
If $\mathcal{M}$ is $k$-curvature homogeneous, then $\mathcal{U}_{m,P}^k:=(T_PM,g_P,R_P,....,\nabla^kR_P)$ is a $k$-model
for $\mathcal{M}$ for any point $P\in M$; conversely, if $\mathcal{M}$ admits a $k$-model, then $\mathcal{M}$ is
$k$-curvature homogeneous.

There are a number of important results in this area in the Riemannian setting ($p=0$). 
Takagi \cite{T74} was the first to exhibit $0$-curvature homogeneous manifolds which are not locally
homogeneous; his examples were non compact. Subsequently, compact examples were exhibited by Ferus, Karcher, and M\"unzer
\cite{FKM81}.
Tomassini
\cite{T97} studied principal fiber bundles with
$1$ dimensional fiber over a $0$-curvature homogeneous base. Other examples may be found in
\cite{KTV92,KTV92a,T88,Va91}. Tsukada \cite{T88x} classified $0$-curvature homogeneous hypersurfaces of dimension $m\ge4$ in
complete and simply connected Riemannian space forms; the case $m=3$ was subsequently treated by Calvaruso, Marinosci, and
Perrone \cite{CMP00}. Kowalski and Pr\"ufer \cite{KP94} exhibited $4$ dimensional algebraic curvature tensors which are not
realizable by any $0$-curvature homogeneous space. 

Scalar invariants can be obtained by using the Weyl calculus to contract indices in pairs in a polynomial expression involving
the curvature and its higher covariant derivatives. For example, the scalar curvature is defined by setting
$$\tau:=\textstyle\sum_{ijkl}g^{ij}g^{kl}R_{iklj}\,.$$
Clearly if $\mathcal{M}$ is locally homogeneous, then all such
scalar invariants are necessarily constant. 

We summarize some important results in this field in the Riemannian setting:
\begin{theorem}\label{thm-1.1}
\ Let $\mathcal{M}$ be a Riemannian manifold. Then:
\begin{enumerate}
\item {\bf (Tricerri and Vanhecke
\cite{TV86})} If $\mathcal{M}$ is modeled
on an irreducible Riemannian symmetric space $\mathcal{N}$, then $\mathcal{M}$ is locally symmetric and hence locally isometric
to $\mathcal{N}$.
\item {\bf (Singer \cite{S60})} There exists an integer $k_m$ so that if $\mathcal{M}$ is a 
complete simply connected manifold of dimension
$m$ which is $k_m$-curvature homogeneous, then $\mathcal{M}$ is homogeneous.
\item {\bf (Pr\"ufer, Tricerri, and Vanhecke
\cite{PTV96})} If all local scalar Weyl invariants up to
order $\frac12m(m-1)$ are constant on a Riemannian manifold $\mathcal{M}$, then $\mathcal{M}$ is locally homogeneous and
$\mathcal{M}$ is determined up to local isometry by these invariants. 
\end{enumerate}
\end{theorem} 

We remark that Cahen et. al. \cite{CLPT90} used a classification result of Berger to show that if $\mathcal{M}$ is a Lorentzian
($p=1$) manifold which is modeled on an irreducible Lorentzian symmetric space, then $\mathcal{M}$ has constant sectional
curvature. Thus Assertion (1) has a natural, and even stronger, extension to the Lorentzian setting.

Singer established the bound $k_m<\frac12m(m-1)$. Bounds of $3m-5$ and $\frac32m-1$ for $k_m$ have been
established by Yamato \cite{Y89} and Gromov \cite{Gr86}. In the low dimensional setting, K. Sekigawa, H. Suga, and L. Vanhecke
\cite{SSV92,SSV95} showed that $k_3=k_4=1$. We refer to the discussion in Boeckx, L. Vanhecke, and O. Kowalski \cite{BVK96} for
further details concerning $k$-curvature homogeneous manifolds in the Riemannian setting. 

Theorem \ref{thm-1.1} (2) extends to the pseudo-Riemannian setting:
\begin{theorem} {\bf(Podesta and Spiro \cite{PS04})}\label{thm-4.1} There exists an integer $k_{p,q}$ so that if $(M,g)$ is a
complete simply connected pseudo-Riemannian manifold of signature $(p,q)$ which is $k_{p,q}$-curvature homogeneous, then $(M,g)$
is homogeneous.
\end{theorem}

B. Opozda \cite{O97} has established an analogue of this result in the affine setting.

In the Lorentzian setting, examples of curvature homogeneous manifolds which are not locally homogeneous were
constructed by Cahen et. al. \cite{CLPT90}. Subsequently, $1$-curvature homogeneous manifolds which are not locally homogeneous
have been constructed by Bueken and Vanhecke
\cite{BV97}; we also refer to related work of Bueken and Djori\'c
\cite{BD00}. These examples are important since they show the results of \cite{SSV92,SSV95} do not extend to the indefinite
setting. Pravda, Pravdov\'a, Coley, and Milson \cite{PPCM02} exhibited Lorentz manifolds all of whose scalar Weyl
invariants vanish and which are not locally homogeneous; thus Theorem \ref{thm-1.1} (3) is false in this setting. 

Not as much is known in the higher signature context. The authors \cite{GS04} exhibited a family
of complete
$1$-curvature homogeneous pseudo-Riemannian manifolds of signature $(r+1,r+1)$ on $\mathbb{R}^{2r+2}$ for
$r\ge2$ which were $0$-modeled on an irreducible symmetric space and which were not $2$-curvature homogeneous (and hence not
homogeneous); two other families of $0$-curvature pseudo-Riemannian manifolds were also exhibited that are $0$-modeled on
irreducible symmetric spaces. Thus Theorem \ref{thm-1.1} (1) fails in the higher signature
setting. We also refer to
\cite{GS03} for other examples of $0$-curvature homogeneous pseudo-Riemannian manifolds.

Let $k=p+2\ge2$ be given. In this paper, we will exhibit a family of complete neutral signature metrics $g_{2p+6,f}$ on
$\mathbb{R}^{2p+6}$ which are $k$-curvature homogeneous but not locally homogeneous for generic values of $f$. We shall be
defining a number of tensors. To simplify the discussion, we shall only give the non-zero entries in these tensors up to the
usual $\mathbb{Z}_2$ symmetries. 

Introduce coordinates $(x,y,z_0,...,z_p,\bar x,\bar y,\bar z_0,...,\bar z_p)$ on
$\mathbb{R}^{2p+6}$. Let
$f=f(y)$ be a smooth function on $\mathbb{R}$. Let $\mathcal{M}_{2p+6,f}:=(\mathbb{R}^{2p+6},g_{2p+6,f})$ be the
pseudo-Riemannian manifold of balanced (i.e. neutral) signature
$(p+3,p+3)$ where:
\begin{eqnarray*}
&&F_{2p+6,f}(y,\vec z):=f(y)+yz_0+y^2z_1+...+y^{p+1}z_p,\\
&&g_{2p+6,f}(\partial_{z_i},\partial_{\bar z_j})=\delta_{ij},\ \ 
  g_{2p+6,f}(\partial_x,\partial_{\bar x})=1,\\
&&g_{2p+6,f}(\partial_y,\partial_{\bar y})=1,\ \ \text{and}\ \ 
g_{2p+6,f}(\partial_x,\partial_x)=-2F_{2p+6,f}(y,\vec z)\,.
\end{eqnarray*}

Choose a basis $\mathcal{B}$ for $\mathbb{R}^{2p+6}$ of the form
 $$\mathcal{B}:=\{X,Y,Z_0,...,Z_p,\bar X,\bar Y,\bar Z_0,...,\bar Z_p\}\,.$$ 
Consider the models $\mathcal{U}_{2p+6}^i:=(\mathbb{R}^{2p+6},g_{2p+6},A_{2p+6}^0,...,A_{2p+6}^i)$ for $0\le i\le p+2$ where the inner product
$g_{2p+6}$ and the tensors $A_{2p+6}^i\in\otimes^{4+i}(\mathbb{R}^{2p+6})^*$ have non-zero components
\begin{equation}\label{eqn-1.a}\begin{array}{l}
g_{2p+6}(X,\bar X)=g_{2p+6}(Y,\bar Y)=g_{2p+6}(Z_i,\bar Z_i)=1,\ \ 
A_{2p+6}^0(X,Y,Z_0,X)=1,\\
A_{2p+6}^1(X,Y,Z_1,X;Y)=A_{2p+6}^1(X,Y,Y,X;Z_1)=1,\vphantom{\vrule height 12pt}\\
A_{2p+6}^2(X,Y,Z_2,X;Y,Y)=A_{2p+6}^2(X,Y,Y,X;Z_2,Y)\vphantom{\vrule height 12pt}\\
\qquad=A_{2p+6}^2(X,Y,Y,X;Y,Z_2)=1,...\vphantom{\vrule height 12pt}\\
A_{2p+6}^p(X,Y,Z_p,X;Y,...,Y)=A_{2p+6}^p(X,Y,Y,X;Z_p,Y,...,Y)\vphantom{\vrule height 12pt}\\
\qquad=...=A_{2p+6}^p(X,Y,Y,X;Y,...,Y,Z_p)=1,\vphantom{\vrule height 12pt}\\
A_{2p+6}^{p+1}(X,Y,Y,X;Y,...,Y)=1,\ \ \text{and}\ \ 
A_{2p+6}^{p+2}(X,Y,Y,X;Y,...,Y)=1\,.
\end{array}\end{equation}

\begin{theorem}\label{thm-1.2}\ 
\begin{enumerate}
\item All geodesics in $\mathcal{M}_{2p+6,f}$ extend for infinite time.
\item $\exp_P:T_P\mathbb{R}^{2p+6}\rightarrow\mathbb{R}^{2p+6}$ is a diffeomorphism for any $P\in\mathbb{R}^{2p+6}$.
\item $\mathcal{U}_{2p+6}^p$ is a $p$-model for $\mathcal{M}_{2p+6,f}$.
\item If $f^{(p+3)}>0$ and $f^{(p+4)}>0$, then $\mathcal{U}_{2p+6}^{p+2}$ is a $p+2$-model for $\mathcal{M}_{2p+6,f}$.
\end{enumerate}
\end{theorem}

It is convenient to work in the affine setting. Let
$$\mathcal{R}(X,Y):=\nabla_X\nabla_Y-\nabla_Y\nabla_X-\nabla_{[X,Y]}$$ 
be the curvature operator defined by a torsion free
connection
$\nabla$ on the tangent bundle of a smooth manifold $M$. Following Opzoda \cite{O96}, we say that
$(M,\nabla)$ is affine {\it $k$-curvature homogeneous} if given any two points $P$ and $Q$ of $M$, there is a linear isomorphism
$\phi:T_PM\rightarrow T_QM$ so that
$\phi^*\nabla^i\mathcal{R}_Q=\nabla^i\mathcal{R}_P$ for $0\le i\le k$. Taking $\nabla$ to be the
Levi-Civita connection of a pseudo-Riemannian metric then yields that any $k$-curvature homogeneous manifold is necessarily
affine
$k$-curvature homogeneous by simply forgetting the requirement that $\phi$ be an isometry; there is no metric present in the
affine setting.

We say that
$(M,\nabla)$ is {\it locally affine homogeneous} if given any points $P$ and $Q$ in $M$, there is a diffeomorphism $\Phi$ from a
neighborhood of $P$ to a neighborhood of $Q$ so that $\Phi(P)=Q$ and so that $\Phi^*\nabla=\nabla$. If $(M,\nabla)$ is
locally affine homogeneous, necessarily $(M,\nabla)$ is affine $k$-curvature homogeneous for any $k$. Examples of $2$-curvature
homogeneous affine manifolds which are not locally affine homogeneous are known; we refer to the discussion in
\cite{GVV99,KO99,KO00,KO04,O96} for this and related results.

We will show that all the scalar Weyl invariants of $\mathcal{M}_{2p+6,f}$ vanish; these manifolds provide additional
examples showing Theorem \ref{thm-1.1} (3) fails in the higher signature setting. To show that $\mathcal{M}_{2p+6,f}$ is not
locally  homogeneous, we must define a suitable invariant. We assume $f^{(p+4)}>0$ and set
$$\alpha_{2p+6,f}:=\frac{f^{(p+3)}f^{(p+5)}}{f^{(p+4)}f^{(p+4)}}\,.$$
Let $\nabla$ be the Levi-Civita connection of $g_{2p+6,f}$. We will show that $\alpha_{2p+6,f}$ is a local affine invariant of
$(\mathbb{R}^{2p+6},\nabla)$; it is not of Weyl type. For generic
$f$, the zero set of the derivative
$\alpha_{2p+6,f}^\prime$ is discrete and hence $\alpha_{2p+6,f}$ is not constant on any open set; thus, for generic $f$, $\mathcal{M}_{2p+6,f}$ is
not locally affine homogeneous and hence not locally homogeneous; furthermore, the scalar Weyl invariants do not determine
$\mathcal{M}_{2p+6,f}$ up to local isometry. 

\begin{theorem}\label{thm-1.3}
Assume that $f^{(p+3)}>0$ and that $f^{(p+4)}>0$. Then:
\begin{enumerate}
\item All scalar Weyl invariants of $\mathcal{M}_{2p+6,f}$ vanish.
\item If $\mathcal{M}_{2p+6,f}$ is affine $p+3$ curvature homogeneous, then $\alpha_{2p+6,f}$ is constant.
\item If $\phi$ is a local diffeomorphism of $\mathcal{M}_{2p+6,f}$ such that $\phi^*\nabla=\nabla$, then we have that
$\phi^*\alpha_{2p+6,f}=\alpha_{2p+6,f}$. 
\item If $\alpha_{2p+6,f}$ is non-constant, then $\mathcal{M}_{2p+6,f}$ is not
locally affine curvature homogeneous.
\end{enumerate}
\end{theorem}

This theorem provides a lower bound for Singer's constant in the neutral setting by showing that if $p\ge0$, then $k_{p+3,p+3}\ge
p+3$ since $f$ can be chosen so $\mathcal{M}_{2p+6,f}$ is $p+2$-curvature homogeneous but not locally $p+3$-affine curvature
homogeneous. By taking suitable product structures and by using the $3$ dimensional Lorentzian examples \cite{BV97} which are
$1$-curvature homogeneous but not locally homogeneous, one may establish the lower bound
$$k_{p,q}\ge\min(p,q)\,.$$
This also establishes a corresponding lower bound in the affine setting for the Opozda constant \cite{O97}.

There are two special cases which are important. Set
$$
\mathcal{M}_{2p+6}^1:=\mathcal{M}_{2p+6,e^y}\quad\text{and}\quad\mathcal{M}_{2p+6}^2:=\mathcal{M}_{2p+6,e^y+e^{2y}}\,.
$$
\begin{theorem}\label{thm-1.4}
\rm\ \begin{enumerate}
\item $\mathcal{M}_{2p+6}^1$ is a homogeneous space.
\item $\mathcal{M}_{2p+6}^2$ is $2p+2$-modeled on $\mathcal{M}_{2p+6}^1$.
\item $\mathcal{M}_{2p+6}^2$ is not locally $p+3$-affine curvature homogeneous.
\end{enumerate}
\end{theorem}

The Jacobi operator is the self-adjoint operator characterized by the property $g(J(X)Y,Z)=R(Y,X,X,Z)$.
One says that $\mathcal{M}$ is {\it nilpotent Osserman} if $0$ is the only  eigenvalue of the Jacobi operator $J(X)$ for any
tangent vector $X$. If $\{e_1,e_2\}$ is an oriented orthonormal basis for a non-degenerate $2$-plane $\pi$, then the
skew-symmetric endomorphism $\mathcal{R}(\pi):=\mathcal{R}(e_1,e_2)$ is independent of the particular basis chosen. One says that
$\mathcal{M}$ is {\it nilpotent Ivanov-Petrova} if
$0$ is the only eigenvalue of $\mathcal{R}(\pi)$ for any such $\pi$. We refer to \cite{GKV02,Gi01} for a further discussion of
these operators in this context.

\begin{theorem}\label{thm-1.5}  $\mathcal{M}_{2p+6,f}$ is Ricci flat,
nilpotent Osserman, and nilpotent Ivanov-Petrova.
\end{theorem}

Theorem \ref{thm-1.1} (1) fails in this setting. We refer to \cite{GS04} for a further discussion of this phenomena and here
content ourselves with showing:

\begin{theorem}\label{thm-1.6} 
Assume that $f^{(3)}>0$ and $f^{(4)}>0$. Then $\mathcal{M}_{6,f}$ is a $6$ dimensional neutral
signature manifold which is $2$-curvature homogeneous, which is
complete, which is modeled on an irreducible neutral signature symmetric space, all of whose local scalar Weyl invariants
vanish identically, and which is not affine $3$-curvature homogeneous for generic $f$.
\end{theorem}

There is a $4$ dimensional example $\mathcal{M}_{4,f}:=(\mathbb{R}^4,g_{4,f})$ where
$$g_{4,f}(\partial_x,\partial_x)=-2f(y)\ \ \text{and}\ \ 
g_{4,f}(\partial_x,\partial_{\bar x})=g_{4,f}(\partial_y,\partial_{\bar y})=1\,.$$
This example is defined, at least in a formal sense, by setting $p=-1$ in the discussion given above. Assume $f^{(2)}>0$ and
$f^{(3)}>0$. Dunn
\cite{D05} showed that
$\mathcal{M}_{4,f}$ is a
$1$-curvature homogeneous  complete manifold which is
$0$-modeled on an irreducible symmetric space and which is not locally homogeneous for generic $f$.

The remainder of this paper is devoted to the proof Theorems \ref{thm-1.2}-\ref{thm-1.6}. In Section \ref{sect-2}, we determine
the Christoffel symbols of the connection $\nabla$ relative to the coordinate frame and establish Assertions (1) and (2) of
Theorem
\ref{thm-1.2}. In Section
\ref{sect-3}, we compute the curvature of the metric $g_{2p+6,f}$;  Theorem
\ref{thm-1.3} (1) and Theorem
\ref{thm-1.5} follow from this computation. In Section \ref{sect-4}, we prove Assertions (3) and (4) of Theorem \ref{thm-1.2}. In
Section
\ref{sect-5}, we complete the proof of Theorem \ref{thm-1.3}; Theorem \ref{thm-1.4} follows as a scholium to these computations.
We conclude the paper in Section
\ref{sect-6} with the proof of Theorem \ref{thm-1.6}.

\section{The geodesics of $\mathcal{M}_{2p+6,f}$}\label{sect-2}
The non-zero Christoffel symbols of the
first and second kinds are given by:
\begin{eqnarray*}
&&g_{2p+6,f}(\nabla_{\partial_x}\partial_y,\partial_x)=g_{2p+6,f}(\nabla_{\partial_y}\partial_x,\partial_x)=-g_{2p+6,f}(\nabla_{\partial_x}\partial_x,\partial_y)
 \\&&\qquad=-\partial_yF_{2p+6,f},\\
&&g_{2p+6,f}(\nabla_{\partial_{z_i}}\partial_x,\partial_x)=g_{2p+6,f}(\nabla_{\partial_x}\partial_{z_i},\partial_x)
 =-g_{2p+6,f}(\nabla_{\partial_x}\partial_x,\partial_{z_i})\\&&\qquad=-y^{i+1},
\end{eqnarray*}
and by
\begin{eqnarray*}
&&\nabla_{\partial_x}\partial_y
   =\nabla_{\partial_y}\partial_x=-(\partial_yF_{2p+6,f})\partial_{\bar x},\\
&&\nabla_{\partial_x}\partial_x=(\partial_yF_{2p+6,f})\partial_{\bar y}+\textstyle\sum_i
y^{i+1}\partial_{\bar z_i},\\
&&\nabla_{\partial_x}\partial_{z_i}=\nabla_{\partial_{z_i}}\partial_x=-y^{i+1}\partial_{\bar x}\,.
\end{eqnarray*}
This exhibits a crucial feature of this metric:
\begin{equation}\label{eqn-2.a}
\nabla\{\partial_x,\partial_y,\partial_{z_i}\}\in\operatorname{Span}\{\partial_{\bar x},\partial_{\bar
y},\partial_{\bar z_i}\},\ \ \text{and}\ \ 
\nabla\{\partial_{\bar x},\partial_{\bar y},\partial_{\bar z_i}\}=\{0\}\,.
\end{equation}

Assertions (1) and (2) of Theorem \ref{thm-1.2} will follow from the following technical Lemma by setting:
$$\begin{array}{rrrrr}
u_1=x,&u_2=y,&u_3=z_0,&...,&u_{p+3}=z_p,\\
u_{p+4}=\bar x,&u_{p+5}=\bar y,&u_{p+6}=\bar z_0,&...,&u_{2p+6}=\bar z_p\,.
\vphantom{\vrule height 12pt}\end{array}$$

\begin{lemma}\label{lem-2.1}
Let $(u_1,...,u_n)$ be coordinates on $\mathbb{R}^n$. Let $g$ be a pseudo-Riemannian metric on
$\mathbb{R}^n$ so that
$\nabla_{\partial_{u_a}}\partial_{u_b}=\sum_{a,b<c}\Gamma_{ab}{}^c(u_1,...,u_{c-1})\partial_{u_c}$.
Then:\begin{enumerate}
\item $(\mathbb{R}^n,g)$ is a complete pseudo-Riemannian manifold.
\item $\exp_P:T_P\mathbb{R}^n\rightarrow\mathbb{R}^n$ is a diffeomorphism for all $P$ in $\mathbb{R}^n$.
\end{enumerate}
\end{lemma}

\begin{proof} We shall adopt the notational convention that the empty sum is $0$. Let $\gamma(t)=(u_1(t),...,u_n(t))$ be a curve
in
$\mathbb{R}^n$;
$\gamma$ is a geodesic if and only
$$
\ddot u_c(t)+\textstyle\sum_{a,b<c}\dot u_a(t)\dot u_b(t)\Gamma_{ab}{}^c(u_1,...,u_{c-1})(t)=0\,.
$$
We solve this system of equations recursively. Let $\gamma(t;\vec u^{\phantom{.}0},\vec u^{\phantom{.}1})$ be defined by
$$
u_c(t)=u_c^0+u_c^1t-\textstyle\int_0^t\int_0^s\textstyle\sum_{a,b<c}\dot u_a(r)\dot
u_b(r)\Gamma_{ab}{}^c(u_1,...,u_{c-1})(r)drds\,.
$$
Then $\gamma(0;\vec u^{\phantom{.}0},\vec u^{\phantom{.}1})(0)=\vec u^{\phantom{.}0}$ while $\dot\gamma(0;\vec
u^{\phantom{.}0},\vec u^{\phantom{.}1})(0)=\vec u^{\phantom{.}1}$. Thus every geodesic arises in this way so all geodesics
extend for infinite time. Furthermore, given $P,Q\in\mathbb{R}^n$, there is a unique geodesic $\gamma=\gamma_{P,Q}$ so that
$\gamma(0)=P$ and
$\gamma(1)=Q$ where
$$
u_c^0=P_c,\ \ u_c^1=Q_c-P_c+\textstyle\int_0^1\int_0^s\textstyle\sum_{a,b<c}\dot
u_a(r)\dot u_b(r)\Gamma_{ab}{}^c(u_1,...,u_{c-1})(r)drds\,.
$$
This shows that $\exp_P$ is a diffeomorphism from $T_P\mathbb{R}^n$ to $\mathbb{R}^n$.
\end{proof}

\section{The curvature of $\mathcal{M}_{2p+6,f}$}\label{sect-3}
In view of Equation (\ref{eqn-2.a}), in computing curvatures and higher covariant derivatives, only derivatives
of highest weight play a role; we never need to consider quadratic terms in Christoffel symbols. Thus
 the non-zero curvatures are:
$$
R_{2p+6,f}(\partial_x,\partial_y,\partial_y,\partial_x)=(\partial_y)^2F_{2p+6,f},\ \ \text{and}\ \ 
R_{2p+6,f}(\partial_x,\partial_y,\partial_{z_i},\partial_x)=(i+1)y^i\,.
$$

\begin{proof}[Proof of Theorem \ref{thm-1.5}] Let $\xi_i$ be arbitrary tangent vectors. Then:
\begin{eqnarray*}
&&\operatorname{Range}\{\mathcal{R}_{2p+6,f}(\xi_1,\xi_2)\}\subset\operatorname{Span}_{C^\infty}\{\partial_{\bar
x},\partial_{\bar y},\partial_{\bar z_0}, ...,\partial_{\bar z_p}\},\quad\text{and}\\
&&\operatorname{Span}\{\partial_{\bar x},\partial_{\bar y},\partial_{\bar z_0},
...,\partial_{\bar z_p}\}\subset\operatorname{Ker}\{\mathcal{R}_{2p+6,f}(\xi_1,\xi_2)\}\,.
\end{eqnarray*}
Thus $\mathcal{R}_{2p+6,f}(\xi_1,\xi_2)\mathcal{R}_{2p+6,f}(\xi_3,\xi_4)=0$ so
$J_{2p+6,f}(\xi)^2=0$ and
$\mathcal{R}_{2p+6,f}(\pi)^2=0$ for any tangent vector $\xi$ and any non-degenerate $2$-plane $\pi$. Consequently,
$J_{2p+6,f}(\xi)$ and
$\mathcal{R}_{2p+6,f}(\pi)$ have only the eigenvalue 0. \end{proof}

Similarly, the non-zero entries in $\nabla^kR$ for any $k\ge0$ are given by:
\begin{eqnarray*}
&&\nabla^kR_{2p+6,f}(\partial_x,\partial_y,\partial_y,\partial_x;\partial_y,...,\partial_y)=(\partial_y)^{k+2}F_{2p+6,f},\\
&&\nabla^kR_{2p+6,f}(\partial_x,\partial_y,\partial_{z_i},\partial_x;\partial_y,...,\partial_y)
=\partial_{z_i}(\partial_y)^{k+1}F_{2p+6,f},
   \ \ \text{and}\\
&&\nabla^kR_{2p+6,f}(\partial_x,\partial_y,\partial_y,\partial_x;\partial_y,...,\partial_{z_i},...,\partial_y)
   =\partial_{z_i}(\partial_y)^{k+1}F_{2p+6,f}\,.
\end{eqnarray*}

\begin{proof}[Proof of Theorem \ref{thm-1.3} (1)] We
may decompose $T\mathbb{R}^{2p+6}=\mathcal{V}\oplus\bar{\mathcal{V}}$ where
\begin{eqnarray*}
&&\mathcal{V}:=\operatorname{Span}\{\partial_x+\ffrac12g_{2p+6,f}(\partial_x,\partial_x)\partial_{\tilde
x},\partial_y,\partial_{z_0},...,\partial_{z_p}\},\quad
\text{and}\\ 
&&\bar{\mathcal{V}}:=\operatorname{Span}\{\partial_{\bar x},\partial_{\bar y},\partial_{\bar
z_0},...,\partial_{\bar z_p}\}\,.
\end{eqnarray*}
Let $\pi_1$ denote projection on the first factor. There are tensors $A^k\in\otimes^{k+4}\mathcal{V}^*$ so that
$\pi_1^*A^k=\nabla^kR$. Since $\mathcal{V}$ is a totally isotropic subspace, this shows all scalar invariants formed using the
Weyl calculus vanish. \end{proof}

\section{A model for $\mathcal{M}_{2p+6,f}$}\label{sect-4}
We can now make a crucial observation. We have
\begin{equation}\label{eqn-4.a}
\nabla^kR_{2p+6,f}(\partial_x,\partial_y,\partial_{z_i},\partial_x;\partial_y,...,\partial_y)=\left\{
\begin{array}{lll}0&\text{if}&i<k,\\
(k+1)!&\text{if}&i=k\,.
\end{array}\right.
\end{equation}
\begin{proof}[Proof of Theorem \ref{thm-1.2} (3,4)] We shall exploit the upper triangular form of Equation (\ref{eqn-4.a}). Let
$a^i(y,\vec z)$ and
$b_i^j(y,\vec z)$ be smooth functions to be chosen presently. Set
$$X=\partial_x-\ffrac12g_{2p+6,f}(\partial_x,\partial_x)\partial_{\bar x},\ \ 
Y=\partial_y+\textstyle\sum_ja^j\partial_{z_j},\ \ \text{and}\ \ 
Z_i=\textstyle\sum_jb_i^j\partial_{z_j}\,.
$$
Assume the matrix $(b_i^j)$ is invertible; let $(\hat b_i^j)$ be the inverse matrix. Set dually
$$\bar X=\partial_{\bar x},\ \ \bar Y=\partial_{\bar y},\ \ \text{and}\ \ 
\bar Z_i=-\textstyle\sum_ja^j\hat b_j^i\partial_{\bar y}+\textstyle\sum_j\hat b_j^i\partial_{\bar z_j}\,.
$$
This is then a hyperbolic basis, i.e. the first relation of Equation (\ref{eqn-1.a}) holds.

We shall assume the matrix $b_i^j$ is triangular:
$$Z_i=\textstyle\sum_{j\le i}b_i^j\partial_{z_j}\,.$$
The relation $\nabla^kR(X,Y,Y,X;Y,...,Y)=0$ for $0\le k\le p$ leads to the equations:
\begin{eqnarray*}
0&=&\nabla^pR(\partial_x,\partial_y,\partial_y,\partial_x;\partial_y,...)
+(p+1)a^pR(\partial_x,\partial_y,\partial_{z_p},\partial_x;\partial_y,...),\\
0&=&\nabla^{p-1}R(\partial_x,\partial_y,\partial_y,\partial_x;\partial_y,...)
+p\textstyle\sum_{p-1\le i\le p}a^iR(\partial_x,\partial_y,\partial_{z_i},\partial_x;\partial_y,...),\quad...\\
0&=&R(\partial_x,\partial_y,\partial_y,\partial_x)
+\textstyle\sum_{0\le i\le p}a^iR(\partial_x,\partial_y,\partial_{z_i},\partial_x)\,.
\end{eqnarray*}
By Equation (\ref{eqn-4.a}), $\nabla^kR(\partial_x,\partial_y,\partial_{z_k},\partial_x,\partial_y,...)\ne0$ and
thus this triangular system of equations determines the coefficients $a^i$ uniquely.

Similarly, the relations $\nabla^kR(X,Y,Z_j,X;Y,...)=\delta_{jk}$ leads to the equations:
\begin{eqnarray*}
1&=&b_p^p\nabla^pR(\partial_x,\partial_y,\partial_{z_p},\partial_x;\partial_y,...),\\
1&=&b_{p-1}^{p-1}\nabla^pR(\partial_x,\partial_y,\partial_{z_{p-1}},\partial_x;\partial_y,...),\\
0&=&\textstyle\sum_{p-1\le i\le p}
   b_p^i\nabla^{p-1}R(\partial_x,\partial_y,\partial_{z_i},\partial_x;\partial_y,...),\quad...,\\
1&=&b_0^0R(\partial_x,\partial_y,\partial_{z_0},\partial_x),\\
0&=&\textstyle\sum_{0\le i\le 1}
   b_1^i\nabla^{p-1}R(\partial_x,\partial_y,\partial_{z_i},\partial_x;\partial_y,...),\\
0&=&\textstyle\sum_{0\le i\le p}
b_p^iR(\partial_x,\partial_y,\partial_{z_i},\partial_x)\,.
\end{eqnarray*}
This system of equations is trianglar. First solve for $b_p^p$, then for $\{b_{p-1}^{p-1},b_p^{p-1}\}$, and finally for
$\{b_0^0,...,b_p^0\}$. Again, the fact that $\nabla^kR(\partial_x,\partial_y,\partial_{z_k},\partial_y;\partial_y,...)\ne0$
is crucial.

If $k>p$, then the only non-zero component of $\nabla^kR$ is given by
$$\nabla^kR_{2p+6,f}(\partial_x,\partial_y,\partial_y,\partial_x;\partial_y...\partial_y)=f^{(k+2)}(y)\,.$$
There is still a bit of freedom left in the choice of basis. Let $\varepsilon_0$ and $\varepsilon_1$ be non-zero functions. We
set
$$\begin{array}{lllll}
X^1=\varepsilon_0X,&Y^1=\varepsilon_1Y,&Z^1_0=\varepsilon_0^{-2}\varepsilon_1^{-1}Z_0,&...,&Z^1_p=
\varepsilon_0^{-2}\varepsilon_1^{-p-1}Z_p,\\
\bar{X^1}=\varepsilon_0^{-1}\bar X,
&\bar{Y^1}=\varepsilon_1^{-1}\bar Y,&
\bar{Z}^1_0=\varepsilon_0^2\varepsilon_1^1\bar Z_0,&...,&
\bar{Z}^1_p=\varepsilon_0^2\varepsilon_1^{p+1}\bar Z_p\,.\vphantom{\vrule height 12pt}
\end{array}$$
The normalizations of Equation (\ref{eqn-1.a}) are preserved for $\{g_{2p+6,f},R,...,\nabla^pR\}$. Also,
\begin{eqnarray*}
\nabla^{p+1}R_{2p+6,f}(X^1,Y^1,Y^1,X^1;Y^1...Y^1)=\varepsilon_0^2\varepsilon_1^{p+3}f^{(p+3)},\\
\nabla^{p+2}R_{2p+6,f}(X^1,Y^1,Y^1,X^1;Y^1...Y^1)=\varepsilon_0^2\varepsilon_1^{p+4}f^{(p+4)}\,.
  \vphantom{\vrule height 12pt}
\end{eqnarray*}
As $f^{(p+3)}>0$ and $f^{(p+4)}>0$, we may set
$$\varepsilon_1:=\frac{f^{(p+3)}}{f^{(p+4)}}\quad\text{and}\quad
  \varepsilon_0:=\frac1{\{\varepsilon_1^{p+3}f^{(p+3)}\}^{\frac12}}\,.
$$
This shows that $\mathcal{U}_{2p+6}^{p+2}$ is a $p+2$ model for $\mathcal{M}_{2p+6,f}$. \end{proof}

\begin{proof}[Proof of Theorem \ref{thm-1.4} (1)] Suppose we set $f(y)=e^y$, $\varepsilon_0=e^{-y/2}$ and
$\varepsilon_1=1$. Then $\nabla^iR_{2p+6,f}(X^1,Y^1,Y^1,X^1;Y^1...Y^1)=1$ for any $i$. Consequently
$\mathcal{M}_{2p+6}^1$ is a simply connected complete $k$-curvature homogeneous manifold for any $k$. Theorem \ref{thm-4.1} now
implies
$\mathcal{M}_{2p+6}^1$ is homogeneous.
\end{proof}

Note that the full strength of Theorem \ref{thm-4.1} is not necessary. Results of Belger and 
Kowalski \cite{BeKo94} show an analytic pseudo-Riemannian manifold which is $k$-curvature homogeneous for all $k$
is locally homogeneous; in our setting the exponential coordinates are  analytic diffeomorphisms so the qualifier `local'
can be removed.

\section{A local invariant}\label{sect-5} Let $k\ge p+1$. Define a generalization of the classical Jacobi operator by setting
$$J_{k,2p+6,f}(Y):X\rightarrow \nabla^k_{Y,...,Y}R_{2p+6,f}(X,Y)Y\,.$$
 Expand $X=a\partial_x+b\partial_y$ and $Y=c\partial_x+d\partial_y$. Then
$$J_{k,2p+6,f}(Y)X=(ad-bc)d^kf^{(k+2)}(d\partial_{\bar x}-c\partial_{\bar y})\,.$$

\begin{proof}[Proof of Theorem \ref{thm-1.3} (2)] Choose $\{Y,X\}$ so $J_{p+1,2p+6,f}(Y)X\ne0$. 
Then necessarily $d\ne0$ and $(ad-bc)\ne0$. 
Let $h$ be {\bf any} Riemannian metric on $\mathcal{M}_{2p+6,f}$;
\begin{eqnarray*}
&&\frac{h(J_{p+1,2p+6,f}(Y)X,J_{p+3,2p+6,f}(Y)X)}{h(J_{p+2,2p+6,f}(Y)X,J_{p+2,2p+6,f}(Y)X)}\\
&=&
\frac{(ad-bc)^2d^{2p+4}f^{(p+3)}f^{(p+5)}}{(ad-bc)^2d^{2p+4}f^{(p+4)}f^{(p+4)}}
\frac{h(d\partial_{\bar x}-c\partial_{\bar y},d\partial_{\bar x}-c\partial_{\bar y})}
{h(d\partial_{\bar x}-c\partial_{\bar y},d\partial_{\bar x}-c\partial_{\bar y})}\\
&=&\alpha_{2p+6,f}\,.
\end{eqnarray*}
Thus $\alpha_{2p+6,f}$ is an affine invariant of
$\{\nabla^{p+1}\mathcal{R},\nabla^{p+2}\mathcal{R},\nabla^{p+3}\mathcal{R}\}$.
\end{proof}

\begin{proof}[Proof of Theorem \ref{thm-1.4} (2,3)] If we set $f=e^y+e^{2y}$, then $\alpha_{2p+6,f}$ is not locally constant so
$\mathcal{M}_{2p+6}^2$ is not locally $p+3$-affine curvature homogeneous. It is, however, $p+2$-curvature modeled on
$\mathcal{M}_{2p+6}^1$.
\end{proof}

\section{Irreducibility}\label{sect-6} We restrict to the case $p=0$. Set $f=0$ to define
$\mathcal{M}_{6,0}$.
The discussion in Section \ref{sect-2} then yields that $\mathcal{M}_{6,0}$ is complete. The computations of Section \ref{sect-3} show
$\nabla R_{g_{6,0}}=0$ so $\mathcal{M}_{6,0}$ is a symmetric space. Furthermore the discussion of Section \ref{sect-4} shows that
$\mathcal{U}_6^0$ is a $0$-model for $\mathcal{M}_{6,0}$. Thus $\mathcal{M}_{6,0}$ is a $0$-model for $\mathcal{M}_{6,f}$. We
complete the proof of Theorem \ref{thm-1.6} by showing that $\mathcal{U}_6^0$ is irreducible as the other assertions then follow.

Let $Z=Z_0$ and $\bar Z=\bar Z_0$. Let $\mathbb{R}^3=\operatorname{Span}\{X,Y,Z\}$. We consider an affine model
$\mathcal{V}=(\mathbb{R}^3,B)$ where
$B\in\otimes^4(\mathbb{R}^3)^*$ is defined by
$$B(X,Y,Z,X)=1\,.$$
\begin{lemma}\label{lem-6.1} The affine model $\mathcal{V}$ is irreducible.
\end{lemma}

\begin{proof} Suppose a non-trivial decomposition $\mathbb{R}^3=V_1\oplus V_2$ induces a corresponding
decomposition $B=B_1\oplus B_2$. Assume the notation chosen so $\dim(V_1)=2$ and $\dim(V_2)=1$. Let $0\ne\xi\in V_2$.
Since $\dim(V_2)=1$, $B_2=0$ so $B(\eta_1,\eta_2,\eta_3,\xi)=0$ for all
$\eta_i\in\mathbb{R}^3$. We expand $\xi=aX+bY+cZ$. We then have
$$a=B(\xi,Y,Z,X)=0,\ \ b=B(X,\xi,Z,X)=0,\ \ \text{and}\ \ c=B(X,Y,\xi,X)=0\,.$$
Thus $\xi=0$ which is false. This contradiction proves the Lemma.\end{proof}

Let $\pi$ be the natural projection from $\mathbb{R}^6$ to $W:=\mathbb{R}^6/\mathcal{K}$ where
$$
\mathcal{K}:=\{\xi\in\mathbb{R}^6:A_6^0(\eta_1,\eta_2,\eta_3,\xi)=0\quad\forall\quad\eta_i\in\mathbb{R}^3\}
=\operatorname{Span}\{\bar X,\bar Y,\bar Z\}\,.
$$
We suppose $\mathcal{U}_6^0$ is reducible and argue for a contradiction. Let
$\mathbb{R}^6=V_1\oplus V_2$ be a non-trivial decomposition with a corresponding decomposition
\begin{equation}\label{eqn-6.a}
g_{6,0}=g_{6,0,1}\oplus g_{6,0,2}\quad\text{and}\quad A_6^0=A_{6,1}^0\oplus A_{6,2}^0\,.
\end{equation}
This also induces a decomposition
$\mathcal{K}=\mathcal{K}_1\oplus\mathcal{K}_2$. We set $W_i:=V_i/\mathcal{K}_i$ to decompose $W=W_1\oplus W_2$ and
$B=B_1\oplus B_2$. By Lemma \ref{lem-6.1}, this decomposition is trivial; we choose the notation so $W_2=\{0\}$ and hence
$V_2\subset\mathcal{K}_2\subset\mathcal{K}$. Since $\mathcal{K}$ is a null subspace,  $g_{6,0,2}$ is trivial. This is a
contradiction as
$g_{6,0}=g_{6,0,1}\oplus g_{6,0,2}$ and $g_{6,0}$ is non-singular. This contradiction completes the
proof of Theorem \ref{thm-1.6}.\hfill\qedbox

\section*{Acknowledgments} Research of P. Gilkey partially supported by the
Max Planck Institute in the Mathematical Sciences (Leipzig). Research of S. Nik\v cevi\'c partially supported by MM 1646
(Srbija). We are grateful to Professors O. Kowalski and L. Vanhecke for introducing us to this area in the first instance and to
Prof. Garc\'{\i}a--R\'{\i}o for suggesting this as a fruitful area of inquiry. We also acknowledge with pleasure helpful
conversations with Dra. S.  L\'opez Ornat.

\end{document}